\documentclass[11pt]{article}
 \usepackage{graphicx,amsmath,amsfonts,amssymb,color,mathrsfs}

\usepackage{theorem}

\def\P{{\mathbb P}}

\catcode`@=12

\newtheorem{thm}{\noindent Theorem}[section]
\newtheorem{lem}{\noindent Lemma}[section]

\theoremheaderfont{\normalfont\bfseries}
\theorembodyfont{\slshape} {\theorembodyfont{\rmfamily}} {\theorembodyfont{\rmfamily}
\newtheorem{defn}{\noindent Definition}[section]}
{\theorembodyfont{\rmfamily} } {\theorembodyfont{\rmfamily}
}
{\theorembodyfont{\rmfamily} } {\theorembodyfont{\rmfamily}
\newtheorem{rem}{\noindent Remark}[section]}
{\theorembodyfont{\rmfamily} } {\theorembodyfont{\rmfamily} } {\theorembodyfont{\rmfamily}
}
 \def\beqlb{\begin{eqnarray}}\def\eeqlb{\end{eqnarray}}
 \def\beqnn{\begin{eqnarray*}}\def\eeqnn{\end{eqnarray*}}
 \numberwithin{equation}{section}
\def\qed{\hfill$\square$\smallskip}

\begin{document}
\title{Exact Tail Asymptotics for a Discrete-Time Preemptive Priority Queue
\footnotetext{E-mail addresses: yy\_0605xx@126.com (Song, Y.),  math\_lzm@csu.edu.cn (Liu, Z.),\;mathdsh@gmail.com (Dai, H.)}
}
\author{\small Yang Song $^1$,
\small Zaiming Liu $^1$,\;Hongshuai Dai$^{2}$ \thanks{Corresponding author}
\\ \small $^1$ School of Mathematics and Statistics, Central South University,  Changsha, 410075 China
\\ \small $^2$ School of Statistics, Shandong University  of Finance and Economics, Jinan, 250014 China
 }
\maketitle
\begin{abstract}
\ In this paper, we consider a discrete-time preemptive priority queue with different service completion probabilities for two classes of customers, one with high-priority and the other with low-priority. This model corresponds to the classical preemptive priority queueing system with two classes of independent Poisson customers and a single exponential server. Due to the possibility of customers' arriving and departing at the same time in a discrete-time queue, the model considered in this paper is more complicated than the continuous-time model. In this model, we focus on the characterization of the exact tail asymptotics for the joint stationary distribution of the queue length of the two types of customers, for the two boundary distributions and for the two marginal distributions, respectively. By using generating functions and the kernel method, we get the exact tail asymptotic properties along the direction of the low-priority queue, as well as along the direction of the high-priority queue.

\end{abstract}

\small {{\bf MSC(2000):}  60K25, 60J10.

{\bf Keywords:} Discrete-time queue $\cdot$ Stationary distribution $\cdot$ Kernel method $\cdot$ Exact tail asymptotics}
\section{Introduction}
Preemptive priority queueing system could be meaningful, sometimes even vital in real life, such as at bank counters or in the emergency room. Miller \cite{M1981} identified a special structure on the rate matrix, and suggested an efficient computational scheme for $M/M/1$ priority queue. Since then, the priority queueing systems have drew a lot of interests, and many researchers have been working hard on different topics of preemptive queues, for example, Gail, Hantler and Taylor \cite{GHT1988, GHT1992}, Kao and Narayanan \cite{KN1990}, Takine \cite{T1996}, Alfa \cite{A1998}, Isotupa and Stanford \cite{IS2002}, Alfa, Liu and He \cite{ALH2003}, Drekic and Woolford \cite{DW2005}, Zhao et al. \cite{Z2006} and so on.

While the stationary distribution is very important for a queueing model to characterize its performance, in most cases, it is difficult to obtain an explicit expression of the stationary distribution. This motivates the study of exact tail asymptotic behavior of the stationary distribution because the property of the exact tail asymptotics often leads to the performance bounds, approximations and other properties of the queueing model. In recent years, many new results on exact tail asymptotics for the continuous-time models have been established, including that Li and Zhao \cite{LZ2009} obtained the tail asymptotic results for the classical continuous-time preemptive priority queueing model. However, in contrast to the extensive studies of tail asymptotics for the continuous-time queueing systems, there has been yet little systematic investigation on the exact tail asymptotics for the discrete-time queueing systems. The main reason is that, in a discrete-time queueing system, all queueing activities (e.g.  arrivals and departures) could occur simultaneously, which results in high complexity and difficulty of the analysis. Considering the importance of the discrete-time queueing systems in both theory and practice, it is interesting and important to study the tail asymptotic behaviors for those discrete-time models. Recently, Xue and Alfa \cite{XA2005} discussed a discrete-time priority $BMAP/PH/1$ queue, and obtained the exact tail asymptotics in the marginal distribution for the low-priority queue. We also note that Xu \cite{X2010} extended the continuous-time model studied in Li and Zhao \cite{LZ2009} to a special kind of discrete-time model, and studied the tail asymptotics for this special model.  Inspired by the above works, in this paper, we will generalize the model in Li and Zhao \cite{LZ2009} to a more general discrete-time model, i.e., the discrete-time preemptive priority queue with single server and two types of customers, and study its tail asymptotic properties. We will carry out the analysis through the Kernel method, the details of which could be found in Fayolle, Iasnogorodski and Malyshev \cite{FIM1999}, Li and Zhao \cite{LZ2009,LZ2011,LZ2012}, Li, Tavakoli and Zhao \cite{LTZ2013}.

The rest of the paper is organized as follows. Section 2 provides the model description and the fundamental form. Discussions about the fundamental form and the kernel equation are presented in Section 3, where the expressions of the generating functions are also obtained. In Section 4, the singularity analysis and the Tauberian-like theorem used to determine the tail asymptotics are provided. Sections 5 and 6 give the details about proofs of exact tail asymptotics along both queue directions.

\section{Model description and fundamental form}

As mentioned, we focus on the discrete-time preemptive priority queue with two classes of customers, one with high-priority and the other with low-priority. Both types of customers arrive independently according to two Bernoulli processes with probabilities $p$ and $q$, and their service times follow geometrical distributions with parameters $\mu_{h}$ and $\mu_{l}$, respectively. The service rule is first in, first out (FIFO) for both classes. With the preemptive rule, the service of a low-priority customer is interrupted upon the arrival of a high-priority customer. The interrupted low-priority customer will stay at the head of the waiting line to restart its service immediately after the last high-priority customer in the system completes its service. All processes are mutually independent. Let $Q_{1}(n)$ and $Q_{2}(n)$ be the number of high- and low-priority customers in the system at the time slot division point $n$, including the one being served, respectively, then we develop a discrete-time Markov chain $\{Q_{1}(n),Q_{2}(n)\}$. Without loss of generality, we assume that $p+q+\mu_{h}+\mu_{l}=1$. Denote $\rho_{h}=\frac{p}{\mu_{h}}$, $\rho_{l}=\frac{q}{\mu_{l}}$ and $\bar{x}=1-x$ for any real number $x\in[0,1]$. One can easily get that the system is stable if $\rho=\rho_{h}+\rho_{l}<1$, which also implies $p<\bar{p}$ and $q<\bar{q}$. Under this condition, for any $i,j=0,1,\cdots$, let $\pi_{ij}$ be the joint stationary distribution of the number of high- and low-priority customers in the system. We also use the following convention: for two functions $f(n)$ and $g(n)$ of nonnegative integers, $f(n)\sim g(n)$ means that $\lim_{n\rightarrow \infty}\frac{f(n)}{g(n)}=1$.

Different from the continuous-time queueing system, all queueing activities such as the customers' potential arrival and departure in a discrete-time queueing system could occur at the same time, so for mathematical clarity, we consider an early arrival system (EAS) (see \cite{HJ1983}) in this paper. Figure 1 explicates the occurrence order of the potential arrival and departure.

\setlength{\unitlength}{0.1in}
\begin{picture}(3,8)
 \put(0,0){\line(1,0){55}}
 \put(0,0){\makebox(5,0.5)[1]{\line(0,1){0.5}}}
 \put(0,-1){\makebox(5,5)[1]{\hbox{$n^{-}$}}}
\put(0,0){\makebox(12,-5)[1]{\vector(0,-1){4}}}
\put(0,0){\makebox(12,0)[1]{\hbox{$\star$}}}
\put(0,0){\makebox(12,-11)[1]{\hbox{$D$}}}

 \put(0,0){\makebox(21,0.5)[1]{\line(0,1){0.5}}}
 \put(0,-1){\makebox(21,5)[1]{\hbox{$n$}}}
  \put(0,0){\makebox(26,5)[1]{\vector(0,-1){4}}}
\put(0,0){\makebox(26,0)[1]{\hbox{$\circ$}}}
\put(0,0){\makebox(26,11)[1]{\hbox{$H$}}}
 \put(0,0){\makebox(33,0)[1]{\hbox{$\diamond$}}}
 \put(0,0){\makebox(33,5)[1]{\vector(0,-1){4}}}
 \put(0,0){\makebox(33,11)[1]{\hbox{$L$}}}

 \put(0,0){\makebox(40,0.5)[1]{\line(0,1){0.5}}}
 \put(0,-1){\makebox(40,5)[1]{\hbox{$n^{+}$}}}
\put(0,0){\makebox(60,0.5)[1]{\line(0,1){0.5}}}
 \put(0,-1){\makebox(60,5)[1]{\hbox{$(n+1)^{-}$}}}
   \put(0,0){\makebox(68,-5)[1]{\vector(0,-1){4}}}
\put(0,0){\makebox(68,0)[1]{\hbox{$\star$}}}
\put(0,0){\makebox(68,-11)[1]{\hbox{$D$}}}

 \put(0,0){\makebox(78,0.5)[1]{\line(0,1){0.5}}}
 \put(0,-1){\makebox(78,5)[1]{\hbox{$n+1$}}}
  \put(0,0){\makebox(85,5)[1]{\vector(0,-1){4}}}
\put(0,0){\makebox(85,0)[1]{\hbox{$\circ$}}}
\put(0,0){\makebox(85,11)[1]{\hbox{$H$}}}
 \put(0,0){\makebox(92,0)[1]{\hbox{$\diamond$}}}
\put(0,0){\makebox(92,0)[1]{\hbox{$\diamond$}}}
 \put(0,0){\makebox(92,5)[1]{\vector(0,-1){4}}}
 \put(0,0){\makebox(92,11)[1]{\hbox{$L$}}}
 \put(0,0){\makebox(100,0.5)[1]{\line(0,1){0.5}}}
 \put(0,-1){\makebox(100,5)[1]{\hbox{${\tiny(n+1)^{+}}$}}}
\put(0,-7.5){$\star \text{ \footnotesize~ Potential departure
}~~~~\circ\text{ \footnotesize~
Potential arrival of the high-priority customer
}$}\put(0,-9.5){$\diamond \text{ \footnotesize~ Potential arrival of the low-priority customer} $} \put(15,-12.5){$ \text{ \footnotesize~
{\bf Figure 1:}\ Various time epochs in an early arrival system} $}
\end{picture}

 \vskip3.9cm

For a discrete-time Markov chain, the transition probabilities and the balance equations can be easily obtained from the transition diagram, see Figure 2.
\begin{center}\includegraphics[width=8cm,height=7cm]{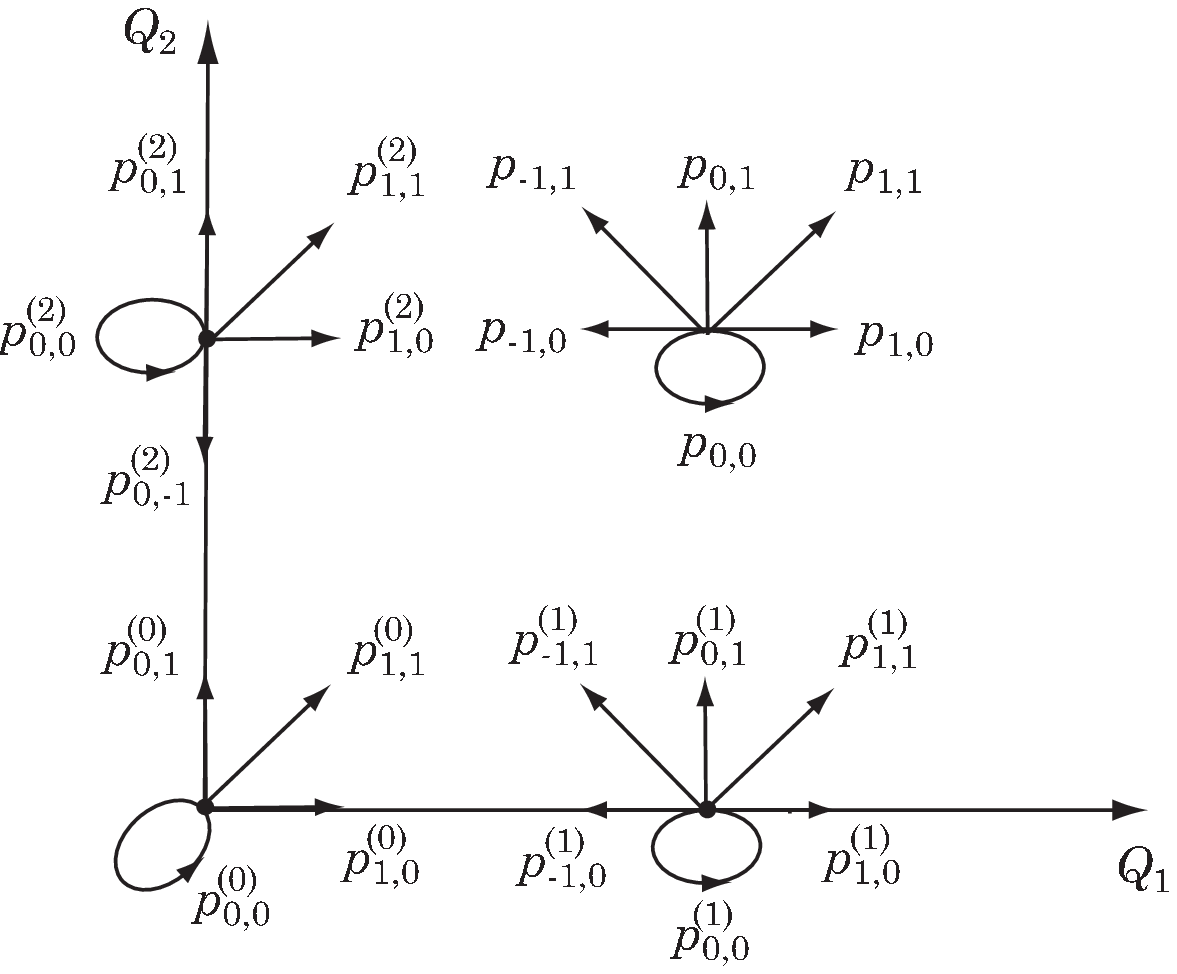}\\
{\footnotesize {\bf Figure 2:} Transition diagram for the discrete-time preemptive priority queue }
\end{center}
The transition probabilities are given by
\beqnn
p_{1,0}&=&p_{1,0}^{(0)}=p_{1,0}^{(1)}=p_{1,0}^{(2)}=p\bar{q}\bar{\mu}_h,
\\p_{1,1}&=&p_{1,1}^{(0)}=p_{1,1}^{(1)}=p_{1,1}^{(2)}=pq\bar{\mu}_h,
\\p_{0,1}&=&p_{0,1}^{(1)}=pq\mu_h+\bar{p}q\bar{\mu}_h,
\\p_{0,1}^{(0)}&=&p_{0,1}^{(2)}=pq\mu_h+\bar{p}q\bar{\mu}_l,
\\p_{-1,1}&=&p_{-1,1}^{(1)}=\bar{p}q\mu_h,
\\p_{-1,0}&=&p_{-1,0}^{(1)}=\bar{p}\bar{q}\mu_h,
\\p_{0,0}&=&p_{0,0}^{(1)}=\bar{p}\bar{q}\bar{\mu}_h+p\bar{q}\mu_h,
\\p_{0,0}^{(0)}&=&\bar{p}\bar{q}+p\bar{q}\mu_h+\bar{p}q\mu_l,
\\p_{0,0}^{(2)}&=&\bar{p}\bar{q}\bar{\mu}_l+p\bar{q}\mu_h+\bar{p}q\mu_l,
\\p_{0,-1}^{(2)}&=&\bar{p}\bar{q}\mu_l.
\eeqnn
Then, the balance equations are
\begin{description}
\item[] \beqlb\label{2-1}
(1-p_{0,0}^{(0)})\pi_{0,0}&=&p_{-1,0}^{(1)}\pi_{1,0}+p_{0,-1}^{(2)}\pi_{01},
\eeqlb
\item[] \beqlb\label{2-1a2}
(1-p_{0,0}^{(1)})\pi_{i,0}&=&p_{1,0}^{(1)}\pi_{i-1,0}+p_{-1,0}^{(1)}\pi_{i+1,0},\quad i\geq1,
\eeqlb
\item[]
\beqlb\label{2-1a3}
(1-p_{0,0}^{(2)})\pi_{0,j}&=&p_{0,1}^{(2)}\pi_{0,j-1}+p_{0,-1}^{(2)}\pi_{0,j+1}+p_{-1,1}\pi_{1,j-1}+p_{-1,0}\pi_{1,j},\quad j\geq1,
\eeqlb
\item[]
\beqlb\label{2-1a4}
(1-p_{0,0})\pi_{i,j}&=&p_{1,0}\pi_{i-1,j}+p_{-1,0}\pi_{i+1,j}+p_{-1,1}\pi_{i+1,j-1}+p_{0,1}\pi_{i,j-1}+p_{1,1}\pi_{i-1,j-1},\quad i\geq1,j\geq1.
\eeqlb
\end{description}
As a specific case of random walks in the quarter plane, we define the following generating functions of the stationary distributions:
\beqnn
\varphi_j(x)&=&\sum_{i=0}^{\infty}\pi_{i,j}x^{i},\quad j\geq0,
\\\psi_i(y)&=&\sum_{j=0}^\infty \pi_{i,j}y^{j},\quad i\geq0,
\\P(x,y)&=&\sum_{i=0}^\infty\sum_{j=0}^\infty \pi_{i,j}x^{i}y^{j}=\sum_{i=0}^{\infty}\psi_i(y)x^i=\sum_{j=0}^{\infty}\varphi_j(x)y^j.
\eeqnn
It is clear that $\varphi_{0}(x)=P(x,0)=P_{1}(x)$ and $\psi_{0}(y)=P(0,y)=P_{2}(y)$. Fayolle, Iasnogorodski and Malyshev \cite{FIM1999} established a functional equation in terms of unknown generating functions, which is often referred as the fundamental form. By using a similar argument, we give the fundamental form as follows, and the bivariate unknown function $P(x,y)$ is connected with two univariate unknown functions $P_{1}(x)$ and $P_{2}(y)$:
\beqlb\label{2-5}
H(x,y)P(x,y)=H_{1}(x,y)P_{1}(x)+H_{2}(x,y)P_{2}(y)+H_{0}(x,y)\pi_{0,0},
\eeqlb
where
\beqnn
H(x,y)&=&-h(x,y),
\\H_1(x,y)&=&-h(x,y)+h_1(x,y)y,
\\H_2(x,y)&=&-h(x,y)+h_2(x,y)x,
\\H_0(x,y)&=&h_0(x,y)xy+h(x,y)-h_1(x,y)y-h_2(x,y)x,
\eeqnn
and
\beqnn
h(x,y)&=&xy\Big(\sum_{i=-1}^1\sum_{j=-1}^1p_{i,j}x^{i}y^{j}-1\Big)=\tilde{a}(y)x^2+\tilde{b}(y)x+\tilde{c}(y),\quad
\\h_1(x,y)&=&x\Big(\sum_{i=-1}^1\sum_{j=0}^1p^{(1)}_{i,j}x^{i}y^{j}-1\Big)=a_1(y)x^2+b_1(y)x+c_1(y),\quad
\\h_2(x,y)&=&y\Big(\sum_{i=0}^1\sum_{j=-1}^1p^{(2)}_{i,j}x^{i}y^{j}-1\Big)=a_2(y)x+b_2(y),\quad
\\h_0(x,y)&=&\sum_{i=0}^1\sum_{j=0}^1p^{(0)}_{i,j}x^{i}y^{j}-1=a_0(y)x+b_0(y).
\eeqnn

After some calculations, we have
\beqlb\label{2-6}
H(x,y)=-y\bigg\{p\bar{\mu}_h(qy+\bar{q})x^{2}+\big[(p\mu_{h}+\bar{p}\bar{\mu}_h)(qy+\bar{q})-1\big]x+\bar{p}\mu_h(qy+\bar{q})\bigg\},
\eeqlb
\beqlb\label{2-7}
H_1(x,y)=0,
\eeqlb
\beqlb\label{2-8}
H_2(x,y)=\bar{p}(qy+\bar{q})\big[(\mu_h-\mu_l)xy-\mu_{h}y+\mu_{l}x\big],
\eeqlb
and
\beqlb\label{2-9}
H_0(x,y)=\bar{p}\bar{q}\mu_l(y-1)x.
\eeqlb
From \eqref{2-6} to \eqref{2-9}, \eqref{2-5} can be rewritten as
\beqlb\label{2-10}
H(x,y)P(x,y)=H_{2}(x,y)\psi_0(y)+H_{0}(x,y)\pi_{0,0}.
\eeqlb

\section{Kernel equation and generating functions}
In this section,  we will consider the {\bf Kernel equation}
\beqnn%\label{a-1}
H(x,y)=0.
\eeqnn

By \eqref{2-6}, we have that
$$H(x,y)=-yK(x,y),$$ where $$K(x,y)=p\bar{\mu}_h(qy+\bar{q})x^{2}
+\big[(p\mu_{h}+\bar{p}\bar{\mu}_h)(qy+\bar{q})-1\big]x+\bar{p}\mu_h(qy+\bar{q}).$$
We call $K(x,y)$   the key kernel. For each fixed $y$, we consider $K(x,y)$ as a quadratic polynomial of $x$ and rewrite it as
\beqnn%\label{3-2}
K(x,y)=a(y)x^{2}+b(y)x+c(y),
\eeqnn
where
\beqnn
a(y)=\frac{\tilde{a}(y)}{y}=p\bar{\mu}_h(qy+\bar{q}),\; b(y)=\frac{\tilde{b}(y)}{y}=(p\mu_h+\bar{p}\bar{\mu}_h)(qy+\bar{q})-1,\;  c(y)=\frac{\tilde{c}(y)}{y}=\bar{p}\mu_h(qy+\bar{q}).
\eeqnn
Let $\Delta(y)$ be the determinant of $K(x,y)=0$, then,
\beqlb\label{3-1}
\nonumber \Delta(y)&=&b(y)^2-4a(y)c(y)
\\&=&(p\mu_h-\bar{p}\bar{\mu}_h)^2(qy+\bar{q})^2-2(p\mu_h+\bar{p}\bar{\mu}_h)(qy+\bar{q})+1.
\eeqlb
Hence, the two solutions to $K(x,y)=0$ are given by
\beqlb\label{3-3}
x_{0}(y)=\frac{1-(p\mu_h+\bar{p}\bar{\mu}_h)(qy+\bar{q})-\sqrt{\Delta(y)}}{2p\bar{\mu}_h(qy+\bar{q})},
\eeqlb
and
\beqlb\label{3-4}
x_{1}(y)=\frac{1-(p\mu_h+\bar{p}\bar{\mu}_h)(qy+\bar{q})+\sqrt{\Delta(y)}}{2p\bar{\mu}_h(qy+\bar{q})}.
\eeqlb
We call $y$ a branch point if it satisfies $\Delta(y)=0$. Here it is easy to get the two branch points
\beqlb\label{3-5}
y_0=\frac{p\mu_h+\bar{p}\bar{\mu}_h-2\sqrt{p\mu_h\bar{p}\bar{\mu}_h}}{(p\mu_h-\bar{p}\bar{\mu}_h)^2q}-\frac{\bar{q}}{q},
\eeqlb
and
\beqlb\label{3-6}
y_1=\frac{p\mu_h+\bar{p}\bar{\mu}_h+2\sqrt{p\mu_h\bar{p}\bar{\mu}_h}}{(p\mu_h-\bar{p}\bar{\mu}_h)^2q}-\frac{\bar{q}}{q}.
\eeqlb
When $y$=0, we have
\beqlb\label{3-7}
x_0=x_0(0)=\frac{r_0}{w},
\eeqlb
\beqlb\label{3-8}
x_1=x_1(0)=\frac{1}{r_0},
\eeqlb
where
\beqlb\label{3-10}
r_0=\frac{1}{x_1(0)},\quad\quad\quad w=\frac{p\bar{\mu}_h}{\bar{p}\mu_h}.
\eeqlb

The following lemma presents the properties of the branch points and the two branches.

\begin{lem}\label{l-1}

\begin{description}
\item[(i)] Suppose that $y_b$ is the root of $b(y)=0$, then $1<y_0<y_b<y_1$.
\item[(ii)] For $-1\leq y\leq1$, we have $x_0(y)<x_1(y)$ and $0<x_0(y)\leq1$. When $y=0$, $0<x_0<1<x_1$.
\end{description}
\end{lem}
{\it Proof:}
(i) One can easily show that $b(y)=(p\mu_h+\bar{p}\bar{\mu}_h)(qy+\bar{q})-1<0$ for $y\leq1$.  This implies $y_b>1$ since that $b(y)$ is a linear function. Moreover, as $a(1)+b(1)+c(1)=1$, we can easily get $$\Delta(1)=[a(1)-c(1)]^2>0.$$  However, $$\Delta(y_b)=-4a(y_b)c(y_b)<0.$$ Noting that $\Delta({\infty})>0$ and $y_0<y_1$, we can conclude $1<y_0<y_b<y_1$.

(ii) It is obvious that $x_0(y)<x_1(y)$. To prove  $x_0(y)>0$ for $y\in[-1,1]$, we need to show $$1-(p\mu_h+\bar{p}\bar{\mu}_h)(qy+\bar{q})>\sqrt{\Delta(y)},$$
which is equivalent to $$4p\mu_h\bar{p}\bar{\mu}_h(qy+\bar{q})^2>0.$$ On the other hand, in order to prove $x_0(y)\leq1$ for $y\in[-1,1]$, it suffices to show $$1-(p\mu_h+\bar{p}\bar{\mu}_h)(qy+\bar{q})-\sqrt{\Delta(y)}\leq2p\bar{\mu}_h(qy+\bar{q})$$ which is equivalent to $$(qy+\bar{q})q(y-1)\leq0.$$
Since $x_1(0)>1$, the proof is completed.
\qed

\begin{rem}
Both the two branches $x_0(y)$ and $x_1(y)$ are analytic in the cut plane $\mathbb{C}_y\setminus[y_0,y_1]$.
\end{rem}

From the balance equations \eqref{2-1} to \eqref{2-1a4},  we can determine the generating functions $\varphi_j(x)$ recursively.

\begin{lem}\label{l-2}
\beqnn%\label{3-11}
\varphi_0(x)=\frac{\pi_{0,0}}{1-r_{0}x},
\eeqnn
and
\beqlb\label{3-12}
\varphi_j(x)=\frac{a_j}{x-x_1}-\frac{q\varphi_{j-1}(x_0)(x+x_0)}{\bar{q}(x-x_1)}-\frac{q(px+\bar{p})(\bar{\mu}_{h}x+\mu_h)}{p\bar{q}\bar{\mu}_h(x-x_1)}\frac{\varphi_{j-1}(x)-\varphi_{j-1}(x_0)}{x-x_0},\quad j=1,2,\ldots
\eeqlb
where
\beqlb\label{3-13}
\pi_{0,0}=\frac{1-\rho}{\bar{p}\bar{q}},
\eeqlb
and
\beqnn
a_j=\frac{\big[\bar{p}\bar{q}(\mu_l-\mu_h)-\bar{p}q\mu_l\big]\pi_{0,j}-\bar{p}\bar{q}\mu_l\pi_{0,j+1}+\bar{p}q(\mu_l-\mu_h)\pi_{0,j-1}-q(\bar{p}\bar{\mu}_h+p\mu_h)\varphi_{j-1}(x_0)}{p\bar{q}\bar{\mu}_h}.
\eeqnn
\end{lem}
{\it Proof:}
It follows from \eqref{2-1a2} that
\beqlb\label{3-14}
\varphi_0(x)=\frac{\bar{p}\bar{q}\mu_h(\pi_{0,0}+\pi_{1,0}x)-\big[1-(p\mu_h+\bar{p}\bar{q})\bar{q}\big]\pi_{0,0}x}{p\bar{q}\bar{\mu}_h(x-x_0)(x-x_1)},
\eeqlb
where
$x_0$ and $x_1$ are given by \eqref{3-7} and \eqref{3-8}, respectively.
By Lemma \ref{l-1} and the fact that $\varphi_0(x)$ is analytic inside the unit circle, we get that $x_0=\frac{r_0}{w}$ is also a zero of the numerator of the function on the right hand side of \eqref{3-14}. Hence
\beqnn
\varphi_0(x)=\frac{\pi_{0,0}}{1-r_{0}x}.
\eeqnn
Now we determine $\pi_{0,0}$. In \eqref{2-10}, let $y=\frac{\mu_{l}x}{\mu_h-(\mu_h-\mu_l)x}$, in which way that the coefficient of $\psi_0(y)$ is zero. Then, $x\rightarrow1$ implies $y \rightarrow1$, and we can obtain $\pi_{0,0}$.

By using the similar argument, we can also get $\varphi_j(x)$ from \eqref{2-1a3} and \eqref{2-1a4}.
\qed

At the end of this section, we determine the generating function $\psi_0(y)$.

\begin{lem}\label{l-3}
Let
\beqlb\label{3-15}
 F(y)=(\bar{p}+\mu_h-2\bar{p}\mu_l)qy^2+\big[(\bar{p}+\mu_{h}-2\bar{p}\mu_l)\bar{q}+2\bar{p}q\mu_{l}-1\big]y+2\bar{p}\bar{q}\mu_l,
\eeqlb
\beqlb\label{3-15-a}
T(y)=F(y)-y\sqrt{\Delta(y)},
\eeqlb
and
 \beqlb\label{3-15-b}
 T^{*}(y)=F(y)+y\sqrt{\Delta(y)},
 \eeqlb
  where $\Delta(y)$ is defined in \eqref{3-1}. Then, for $y\in[-1,1]$, we have
\beqlb\label{3-16}
\psi_0(y)=a\frac{T^*(y)}{(qy+\bar{q})(1-\eta_1{y})}+b\frac{T^*(y)}{(qy+\bar{q})(1-\eta_2{y})},
\eeqlb
where
\beqlb\label{3-17}
a=\frac{\pi_{0,0}}{2\bar{p}\mu_l}\frac{\eta_1}{\eta_1-\eta_2},\quad b=\frac{\pi_{0,0}}{2\bar{p}\mu_l}\frac{\eta_2}{\eta_2-\eta_1},\quad
\eeqlb
\beqlb\label{3-18}
\eta_1=\frac{(1-\mu_h\bar{q}-\bar{p}\bar{q}\bar{\mu}_l-\bar{p}q\mu_l)+\sqrt{(1-\mu_h\bar{q}-\bar{p}\bar{q}\bar{\mu}_l-\bar{p}q\mu_l)^2+4\bar{p}\bar{q}(\mu_h-\bar{p}\mu_l)\bar{\mu}_{l}q}}{2\bar{p}\bar{q}\mu_l},
\eeqlb
\beqlb\label{3-19}
\eta_2=\frac{(1-\mu_h\bar{q}-\bar{p}\bar{q}\bar{\mu}_l-\bar{p}q\mu_l)-\sqrt{(1-\mu_h\bar{q}-\bar{p}\bar{q}\bar{\mu}_l-\bar{p}q\mu_l)^2+4\bar{p}\bar{q}(\mu_h-\bar{p}\mu_l)\bar{\mu}_{l}q}}{2\bar{p}\bar{q}\mu_l}.
\eeqlb
\end{lem}
{\it Proof:}
From \eqref{2-10}, we have
\beqnn
P(x,y)&=&\frac{H_2(x,y)\psi_0(y)+H_0(x,y)\pi_{0,0}}{-yK(x,y)}
\\&=&\frac{\bar{p}(qy+\bar{q})\big[(\mu_h-\mu_l)xy-\mu_{h}y+\mu_{l}x\big]\psi_0(y)+\bar{p}\bar{q}\mu_l(y-1)x\pi_{0,0}}{-y(x-x_0(y))(x-x_1(y))},
\eeqnn
where $x_0(y)$ and $x_1(y)$  are given by \eqref{3-3} and \eqref{3-4}, respectively.  Since $P(x_0(y),y)$ is analytic and nonzero for $-1\leq y\leq1$,  we have
\beqlb\label{3-19-a}
H_2(x_0(y),y)\psi_0(y)+H_0(x_0(y),y)\pi_{0,0}=0.
\eeqlb
Therefore, the  equations \eqref{3-3} and  \eqref{3-19-a} lead to
\beqlb\label{3-20}
\nonumber \psi_0(y)&=&\frac{\bar{q}\mu_l(1-y)x_0(y)\pi_{0,0}}{(qy+\bar{q}){\big\{[(\mu_h-\mu_l)y+\mu_{l}]x_0(y)-\mu_{h}y}\big\}}
\\&=&\frac{2\bar{p}\bar{q}\mu_l(1-y)\pi_{0,0}T^*(y)}{T(y)T^*(y)},
\eeqlb
since  $x_0(y)x_1(y)=1/w$. It is easy to get that $T(1)=0$ and $T^*(-\frac{\bar{q}}{q})=0$.

Next, we  will identify other zeros of $T(y)T^*(y)$.  In fact, we have
\beqnn
T(y)T^*(y)&=&F(y)^2-y^2\Delta(y)
\\&=&4\bar{p}(qy+\bar{q})(y-1)\big[(\mu_{h}-\bar{p}\mu_{l})\bar{\mu}_{l}qy^2+\mu_l(1-\mu_h\bar{q}-\bar{p}\bar{q}\bar{\mu}_l-\bar{p}q\mu_{l})y-\bar{p}\bar{q}\mu_{l}^2\big]
\\&=&4\bar{p}(qy+\bar{q})(y-1)f(y),
\eeqnn
where
\beqlb\label{3-21}
\nonumber f(y)&=&(\mu_{h}-\bar{p}\mu_{l})\bar{\mu}_{l}qy^2+\mu_l(1-\mu_h\bar{q}-\bar{p}\bar{q}\bar{\mu}_l-\bar{p}q\mu_{l})y-\bar{p}\bar{q}\mu_{l}^2
\\&=&-\bar{p}\bar{q}\mu_{l}^2(1-\eta_1y)(1-\eta_2y).
\eeqlb
Obviously, we have $\eta_2<0<\eta_1$, and $\eta_1$, $\eta_2$ are two non-unit zeros of the denominator of the generating function $\psi_0(y)$. From Lemma 3.1 we know that $1<y_0<y_1$, so $T^*(y)$ is analytic in $[-1,1]$. By substituting
\beqnn
T(y)T^{*}(y)=4\bar{p}^{2}\bar{q}\mu_{l}^2(qy+\bar{q})(1-y)(1-\eta_{1}y)(1-\eta_{2}y)
\eeqnn
into $\psi_0(y)$, and then using partial fractions, we can get that \eqref{3-16} holds for $-1\leq y\leq1$, which completes the proof of the lemma.\qed

Next, we will study the exact tail asymptotic behaviors of this queueing system and only consider the case that $\mu_l\leq\mu_h$. In fact, it is not easy  to characterize the exact properties of the singularities in the case of $\mu_l>\mu_h$. It is much complicated and depends on the values of the parameters.

\section{Analysis of singularities}
The analysis of the exact tail asymptotics along the low-priority queue direction in the stationary distribution $\pi_{i,j}$ relies on the analysis of the singularities of the generating function $\psi_0(y)$, which is the focus of this section. Then according to the detailed asymptotic property at the dominant singularity of the generating function, asymptotics of the coefficients of the generating functions will be obtained by using the Tauberian-like theorem. Before giving the Tauberian-like theorem, we first introduce the definition of $\Delta-$ domain in \cite{FS2009}.

\begin{defn}
For given numbers $\epsilon>0$ and $\phi$ with $0<\phi<\frac{\pi}{2}$, the open domain $\Delta(\phi,\epsilon)$ is defined by
\beqnn
\Delta(\phi,\epsilon)=\big\{z\in\mathbb{C}:\,|z|<1+\epsilon, z\neq 1,\,|Arg(z-1)|>\phi\big\}.
\eeqnn
A domain is a $\Delta-$ domain at 1 if it is a $\Delta(\phi,\epsilon)$ for some $\epsilon>0$ and $0<\phi<\frac{\pi}{2}$.
 For a complex number $\xi\neq0$, a $\Delta-$ domain at $\xi$ is defined as the image $\xi\cdot\Delta(\phi,\epsilon)$ of a
 $\Delta-$ domain $\Delta(\phi,\epsilon)$ at 1 under the mapping $z\to\xi z$. A function is called $\Delta-$ analytic if it is analytic in some
 $\Delta-$ domain.
\end{defn}
\begin{rem}\label{a-rem2}
When we use the Tauberian-like theorem to obtain the exact asymptotics,
we often need to consider the region $\Delta(\phi,\epsilon)$. The region $\Delta(\phi,\epsilon)$
is a dented disk. In the sequel, if not otherwise stated,
the limit of a $\Delta-$ analytic function is always taken in $\Delta(\phi,\epsilon)$.
\end{rem}

The following Tauberian-like theorem is from Bender \cite{B1974}, which can also be found in Flajolet and Sedgewick \cite{FS2009}, Li and Zhao \cite{LZ2012}.
\begin{thm}\label{4-thm1}
{\bf (Tauberian-like theorem for single singularity)}
Let $A(z)=\sum_{n\geq 0}a_n z^n$ be analytic at zero with the radius of convergence $R$.
Suppose that $R$ is a singularity of $A(z)$ on the circle of convergence such that $A(z)$ can be continued to a $\Delta-$ domain at $R$. If for a real number
$\alpha\notin \{0,-1,-2,...\}$,

\beqnn
\lim_{z\to R}(1-z/R)^{\alpha}A(z)=g,
\eeqnn
where $g$ is a  non-zero constant. Then
\beqnn
a_n \sim \frac{g}{\Gamma(\alpha)}n^{\alpha-1}R^{-n},
\eeqnn
where $\Gamma(\alpha)$ is the value of Gamma function at $\alpha$.
\end{thm}

To apply the above Tauberian-like theorem, we only need to pay attention to the singularities with modulus greater than 1. For the generating function $\psi_0(y)$, the following Key Lemma identifies the cases when the dominant singularity is a pole or not a pole.
\begin{lem}\label{lem4-1}
\begin{description}
\item[(i)] If $F(y_0)\neq0$, then $1<1/\eta_1<y_0$.
\item [(ii)] If $\rho=\rho_h+\rho_l<1$, then $T^{'}(1)<0$.
\end{description}
\end{lem}
{\it Proof:}
(i) It is easy to get that $f(0)=-\bar{p}\bar{q}\mu_l^2<0$ and $f(1)=\mu_l\mu_h(\rho-1)<0$. From \eqref{3-18}, \eqref{3-19} and \eqref{3-21}, we can also get that $1/\eta_2<0$ and $1/\eta_1>1$. If $F(y_0)\neq0$, then, from $T(y_0)T^*(y_0)=F(y_0)^2>0$, we get $f(y_0)>0$. Since the branch point $y_0>1$,  $1<1/\eta_1<y_0$.

(ii) Noting that
\beqnn
\frac{dT(y)}{dy}=\frac{d F(y)}{dy}-\sqrt{\Delta(y)}-\frac{y}{2\sqrt{\Delta{(y)}}}\frac{d\Delta(y)}{dy},
\eeqnn
we get that $T^{'}(1)=-2\overline{p}\mu_l(1-\rho)/(1-\rho_h)<0$, since $\rho=\rho_h+\rho_l<1$.
\qed

\begin{lem}\label{KL}
{\bf (Key Lemma)}
For the property of $1/\eta_1$, there are three cases:
\begin{description}
\item[(i)]If $F(y_0)>0$, then $1<1/\eta_1<y_0$ and $1/\eta_1$ is a zero of $T(y)$, but not $T^*(y)$. Hence $1/\eta_1$ is the dominant singularity of $\psi_0(y)$, which is a pole.
\item[(ii)] If $F(y_0)=0$, then $1<1/\eta_1=y_0$ and $1/\eta_1=y_0$ is a zero of both $T(y)$ and $T^*(y)$. Hence $1/\eta_1=y_0$ is the dominant singularity of $\psi_0(y)$, which is a branch point.
\item[(iii)] If $F(y_0)<0$, then $1<1/\eta_1<y_0$ and $1/\eta_1$ is a zero of $T^*(y)$, but not $T(y)$. Hence $y_0$ is the dominant singularity of $\psi_0(y)$, which is a branch point.
\end{description}
\end{lem}
{\it Proof:}
 We prove the lemma based on the property of $F(y)$. Since $\bar{p}+\mu_h-2\bar{p}\mu_l>0$, $F(0)=2\bar{p}\bar{q}\mu_l>0$ and $F(1)=\mu_h-p>0$, there are only four cases for the quadratic function $F(y)$:
 \\(a) If the roots of $F(y)=0$ are in $(0,1)$, then $T(y)>0$ for $y\in(-\infty,0)$ implies $T^*(\frac{1}{\eta_2})=0$. But from Lemma \ref{lem4-1} we know $1<1/\eta_1<y_0$, so $T^*(\frac{1}{\eta_1})>0$, which would yield $T(\frac{1}{\eta_1})=0$.
 \\(b) If the roots of $F(y)=0$ are in $(1, +\infty)$, then $T(y)>0$ for $y\in(-\infty,0)$ implies $T^*(\frac{1}{\eta_2})=0$. If $F(y_0)>0$, from Lemma \ref{lem4-1} we have $1<1/\eta_1<y_0$, which would yield  $T(\frac{1}{\eta_1})=0$ since $T(1)=0$, $T^{'}(1)<0$ and $T(y_0)>0$. If $F(y_0)=0$, it implies $f(y_0)=0$, then $1<y_0=1/\eta_1$ and $T(\frac{1}{\eta_1})=T^*({\frac{1}{\eta_1}})=0$. If $F(y_0)<0$, then $T^*(y_0)<0$. From $1<1/\eta_1<y_0$ and $T^*(1)>0$, we can easily show that $T^*(\frac{1}{\eta_1})=0$.
 \\(c) If the roots of $F(y)=0$ are in $(-\infty,0)$, then $T^*(y)>0$ for $y\in(0,y_0)$. From Lemma \ref{lem4-1}, we have $T^*(\frac{1}{\eta_1})>0$. Therefore, $T(\frac{1}{\eta_1})=0$. Moreover, from $F(-\infty)>0$, we get that $T(-\infty)>0$. Finally, we also notice that $T^*(-\frac{q}{\bar{q}})=0$, $T(1)=0$ and $T^{'}(1)<0$. Therefore,  we obtain $T(\frac{1}{\eta_2})\neq0$ but $T^*(\frac{1}{\eta_2})=0$.
\\(d) If $F(y)$ has no zeros, then for any $y$ we have $F(y)>0$. In this case, the equations \eqref{3-15-a} and \eqref{3-15-b} yield $T(\frac{1}{\eta_1})=0$ and $T^*(\frac{1}{\eta_2})=0$ directly.

\qed
\begin{rem}
In our model, it is obvious  that $\psi_0(y)$ must be analytic on the whole complex plane except $[y_0,y_1]\cup\{\frac{1}{\eta_1}\}$, so the Tauberian-like theorem will be applied directly in the next section.
\end{rem}
\begin{rem}
In the next section, we will show that the three cases in the Key Lemma  correspond to the three types of exact tail asymptotics along the low priority queue direction: (i) exact geometric; (ii) geometric with a prefactor $j^{-\frac{1}{2}}$; (iii) geometric with a prefactor $j^{-\frac{3}{2}}$.
\end{rem}

\section{Exact tail asymptotics for the low-priority queue}

In this section, we provide the exact tail asymptotic properties in the joint stationary distribution as well as in the marginal distribution. Firstly, we study the tail asymptotics for the boundary probabilities $\pi_{0,j}$,   and then  obtain the exact tail asymptotic characterization in the joint probabilities $\pi_{i,j}$ for any fixed $i\geq1$.

\subsection{Exact tail asymptotics of the boundary probabilities}

Stationary probabilities $\pi_{0,j}$ for $ j\geq 0$ are referred as the boundary probabilities.
In this subsection, we apply properties obtained in the previous section to characterize the asymptotic behavior of $P_2(y)$, i.e., $\psi_0(y)$.
The exact tail asymptotics of boundary probabilities is a direct consequence of the Tauberian-like theorem and the asymptotic behavior of $\psi_0(y)$.

Define
\beqlb\label{5-1}
C_{l,1}&=&\frac{2aF(\frac{1}{\eta_1})}{q\eta_1^{-1}+\bar{q}},
\eeqlb
\beqlb\label{5-2}
C_{l,2}=\frac{a(\bar{p}-\mu_h)qy_0\sqrt{y_0(y_1-y_0)}}{(qy_0+\bar{q})\sqrt{\pi}},
\eeqlb
\beqlb\label{5-3}
C_{l,3}=\bigg[\frac{a}{(qy_0+\bar{q})(1-\eta_{1}y_0)}+\frac{b}{(qy_0+\bar{q})(1-\eta_{2}y_0)}\bigg]\frac{(\mu_h-\bar{p})q}{2\sqrt{\pi}}\sqrt{y_0(y_1-y_0)},
\eeqlb
where
$a$, $b$, $\eta_1$, $\eta_2$, $y_0$, $y_1$ and $F(\cdot)$ are given by \eqref{3-17}, \eqref{3-18}, \eqref{3-19}, \eqref{3-5}, \eqref{3-6}, \eqref{3-15}, respectively.

Now we state the main result of this section.

\begin{thm}\label{5-thm1-a}
For the discrete-time preemptive priority queue with two classes of customers satisfying $\rho<1$, characterizations of the exact tail asymptotics in the boundary stationary distribution along the low-priority queue direction are given as follows, for $i=0$, where $i$ is the number of high-priority customers:
\\(i) (Exact geometric decay) In the region defined by $F(y_0)>0$,
\beqnn
\pi_{0,j}\sim C_{l,1}\eta_{1}^j.
\eeqnn
\\(ii) (Geometric with a prefactor $j^{-\frac{1}{2}}$) In the region defined by $F(y_0)=0$,
\beqnn
\pi_{0,j}\sim C_{l,2}j^{-\frac{1}{2}}y_0^j.
\eeqnn
\\(iii) (Geometric with a prefactor $j^{-\frac{3}{2}}$) In the region defined by $F(y_0)<0$,
\beqnn
\pi_{0,j}\sim C_{l,3}j^{-\frac{3}{2}}y_0^j.
\eeqnn
\end{thm}
{\it Proof:}
 (i) In this case, $1<1/\eta_1<y_0$.  So, it follows from the Key Lemma \ref{KL} that $T^*(-\frac{\bar{q}}{q})=0$, $T(\frac{1}{\eta_1})=0$ and $T^*({\frac{1}{\eta_2}})=0$. Clearly, $\frac{T^*(y)}{(qy+\bar{q})(1-y\eta_1)}$ is analytic in $\Delta(\phi,\epsilon)=\big\{y\eta_1:\,|y\eta_1|<1+\epsilon, y\eta_1\neq 1,\,|Arg(y\eta_1-1)|>\phi,\epsilon>0,0<\phi<\pi/2,\big\}$.

 On the other hand, we have
\beqnn
\displaystyle\lim_{y\eta_{1}\rightarrow 1}(1-y\eta_1)\psi_0(y)=\displaystyle\lim_{y\eta_{1}\rightarrow 1}\frac{aT^*(y)}{qy+\bar{q}}=\frac{2aF(\frac{1}{\eta_1})}{q\eta_1^{-1}+\bar{q}}=C_{l,1},
\eeqnn
where $C_{l,1}$ is given by \eqref{5-1}.

Therefore,  by Theorem \ref{4-thm1},
\beqnn
\pi_{0,j}\sim C_{l,1}\eta_{1}^j.
\eeqnn

(ii) In this case, $1<1/\eta_1=y_0$.  According to the Key Lemma \ref{KL}, $T(\frac{1}{\eta_1})=T^*(\frac{1}{\eta_1})=0$ and $T^{*}(\frac{1}{\eta_2})=0$. Then
\beqnn
\displaystyle\lim_{y\eta_{1}\rightarrow 1}\sqrt{1-y\eta_1}\psi_0(y)&=&\displaystyle\lim_{y\eta_{1}\rightarrow 1}\bigg[a\frac{T^*(y)}{(qy+\bar{q})\sqrt{1-y\eta_1}}+b\frac{T^*(y)\sqrt{1-y\eta_1}}{(qy+\bar{q})(1-y\eta_2)}\bigg]
\\&=&\displaystyle\lim_{y\eta_{1}\rightarrow 1}a\frac{T^*(y)}{(qy+\bar{q})\sqrt{1-y\eta_1}}
\\&=&\displaystyle\lim_{y\eta_{1}\rightarrow 1}\frac{a}{qy+\bar{q}}\frac{F(y)+y\sqrt{\Delta(y)}}{\sqrt{1-y\eta_1}}
\\&=&\displaystyle\lim_{y\eta_{1}\rightarrow 1}\frac{a}{qy+\bar{q}}\bigg[\sqrt{\frac{F(y)}{1-y\eta_1}}\sqrt{F(y)}+y\sqrt{\frac{y_0\Delta(y)}{y_0-y}}\bigg].
\eeqnn
Rewrite $\Delta(y)$  as follows
\beqnn
\Delta(y)=\big[(p\mu_h-\bar{p}\bar{\mu}_h)q\big]^2(y-y_0)(y-y_1).
\eeqnn
Since $F(\frac{1}{\eta_1})=0$, $\frac{F(y)}{1-\eta_1}$ is a polynomial of degree 1. Therefore,
$$\displaystyle\lim_{y\eta_{1}\rightarrow 1}\sqrt{\frac{F(y)}{1-y\eta_1}}\sqrt{F(y)}=0,$$
which implies
\beqnn
\displaystyle\lim_{y\eta_{1}\rightarrow 1}\sqrt{1-y\eta_1}\psi_0(y)=\frac{a(\bar{p}-\mu_h)q y_0}{qy_0+\bar{q}}\sqrt{y_0(y_1-y_0)}.
\eeqnn
By Theorem \ref{4-thm1}, we have
\beqnn
\pi_{0,j}\sim C_{l,2}j^{-\frac{1}{2}}y_0^{-j},
\eeqnn
where $C_{l,2}$ is given by \eqref{5-2}.

(iii) In this case, $1<1/\eta_1<y_0$ and $T^*(\frac{1}{\eta_1})=0$. We have
\beqnn
\frac{d\psi_0(y)}{dy}=T^*{(y)}\frac{d}{dy}\bigg[\frac{a}{(qy+\bar{q})(1-y\eta_1)}+\frac{b}{(qy+\bar{q})(1-y\eta_2)}\bigg]+\bigg[\frac{a}{(qy+\bar{q})(1-y\eta_1)}+\frac{b}{(qy+\bar{q})(1-y\eta_2)}\bigg]\frac{d T^{*}(y)}{dy},
\eeqnn
where
\beqnn
\frac{dT^{*}(y)}{dy}=\frac{d F(y)}{dy}+\sqrt{\Delta(y)}+\frac{y}{2\sqrt{\Delta(y)}}\frac{d\Delta(y)}{dy},
\eeqnn
\beqnn
\frac{d\Delta(y)}{dy}=\big[(p\mu_h-\bar{p}\bar{\mu}_h)q\big]^2(y-y_1+y-y_0).
\eeqnn
Therefore,
\beqnn
\displaystyle\lim_{y\rightarrow y_0}\sqrt{1-\frac{y}{y_0}}\frac{d\psi_0(y)}{dy}&=&\displaystyle\lim_{y\rightarrow y_0}\bigg[\frac{a}{(qy+\bar{q})(1-y\eta_1)}+\frac{b}{(qy+\bar{q})(1-y\eta_2)}\bigg]\sqrt{1-\frac{y}{y_0}}\frac{dT^{*}(y)}{dy}
\\&=&\displaystyle\lim_{y\rightarrow y_0}\bigg[\frac{a}{(qy+\bar{q})(1-y\eta_1)}+\frac{b}{(qy+\bar{q})(1-y\eta_2)}\bigg]\sqrt{1-\frac{y}{y_0}}\frac{y}{2\sqrt{\Delta(y)}}\frac{d\Delta(y)}{dy}
\\&=&\displaystyle\lim_{y\rightarrow y_0}\bigg[\frac{a}{(qy+\bar{q})(1-y\eta_1)}+\frac{b}{(qy+\bar{q})(1-y\eta_2)}\bigg]\frac{y}{2}\sqrt{\frac{y_0-y}{y_0\Delta(y)}}\big[(p\mu_h-\bar{p}\bar{\mu}_h)q\big]^2(y-y_1+y-y_0)
\\&=&\bigg[\frac{a}{(qy_0+\bar{q})(1-y_0\eta_1)}+\frac{b}{(qy_0+\bar{q})(1-y_0\eta_2)}\bigg]\frac{(\mu_h-\bar{p})q}{2}\sqrt{y_0(y_1-y_0)}.
\eeqnn
Finally,
\beqnn
\pi_{0,j}\sim C_{l,3}j^{-\frac{3}{2}}y_0^{-j},
\eeqnn
where $ C_{l,3}$ is given by \eqref{5-3}.
\qed

\subsection{Exact tail asymptotics of the joint probabilities}

In the previous subsection, exact tail asymptotic properties of the boundary probabilities are obtained. In this subsection, we provide the details for the exact tail asymptotic characterizations in the joint probabilities $\pi_{i,j}$ for any fixed $i\geq 1$.

From the  balance equations \eqref{2-1} to \eqref{2-1a4}, we obtain
\beqlb\label{5-7}
b_2(y)\psi_0(y)+y c(y)\psi_1(y)=\bar{p}\bar{q}\mu_l(1-y)\pi_{0,0},
\eeqlb
and
\beqlb\label{5-8}
a(y)\psi_{i-1}(y)+b(y)\psi_i(x)+c(y)\psi_{i+1}(y)=0,\quad i\geq1.
\eeqlb

Now, we are ready to give the main result of this subsection, which shows the exact tail asymptotics for the joint probabilities.
\begin{thm}\label{5-thm1}
For the discrete-time preemptive priority queue with two classes of customers satisfying $\rho<1$, characterizations of the exact tail asymptotics in the joint stationary distribution along the low-priority queue direction are given as follows: for a fixed number $i$ of the high-priority customers, where $i\geq1$:
\\(i) (Exact geometric decay) In the region defined by $F(y_0)>0$,
\beqlb\label{5-9}
\pi_{i,j}\sim C_{i,1}\eta_{1}^j,
\eeqlb
where
\beqnn
C_{i,1}=\frac{1-(p\mu_h+\bar{p}\bar{\mu}_l)(q\eta_1^{-1}+\bar{q})-\bar{p}\mu_l(q+\bar{q}\eta_1)}{\bar{p}\mu_h(q\eta_1^{-1}+\bar{q})}\bigg[\frac{p\bar{\mu}_h}{\bar{p}\mu_h-\bar{p}\mu_{l}(1-\eta_1)}\bigg]^{i-1}C_{l,1}.
\eeqnn
\\(ii) (Geometric with a prefactor $j^{-\frac{1}{2}}$) In the region defined by $F(y_0)=0$,
\beqnn%\label{5-10}
\pi_{i,j}\sim C_{i,2}j^{-\frac{1}{2}}y_0^j,
\eeqnn
where
\beqnn
C_{i,2}=\frac{1-(p\mu_h+\bar{p}\bar{\mu}_l)(q y_0+\bar{q})-\bar{p}\mu_l(q+\bar{q}y_0^{-1})}{\bar{p}\mu_h(q y_0+\bar{q})}\bigg[\frac{p\bar{\mu}_h}{\bar{p}\mu_h-\bar{p}\mu_{l}(1-y_0^{-1})}\bigg]^{i-1}C_{l,2}.
\eeqnn
\\(iii) (Geometric with a prefactor $j^{-\frac{3}{2}}$) In the region defined by $F(y_0)<0$,
\beqnn%\label{5-11}
\pi_{i,j}\sim C_{i,3}j^{-\frac{3}{2}}y_0^j,
\eeqnn
where
\beqnn
C_{i,3}=-\bigg[\frac{b_2(y_0)}{y_{0}c(y_0)}+\frac{h_2(x_0(y_0),y_0)}{y_{0}c(y_0)}(i-1)\bigg]\bigg[\frac{2p\bar{\mu}_h(q y_0+\bar{q})}{1-(p\mu_h+\bar{p}\bar{\mu}_h)(q y_0+\bar{q})}\bigg]^{i-1}C_{l,3}.
\eeqnn
\end{thm}
{\it Proof:}
(i) Noting that $\lim_{y\eta_{1}\rightarrow 1}(1-y\eta_1)\psi_0(y)=C_{l,1}$, by the equations \eqref{5-7} and \eqref{5-8},  we assume that for  $i\geq1$, $$\lim_{y\eta_{1}\rightarrow 1}(1-y\eta_1)\psi_i(y)=C_{i,1}.$$ Then
\beqnn
b_2(\eta_1^{-1})C_{l,1}+\eta_1^{-1}c(\eta_1^{-1})C_{1,1}=0,
\eeqnn
and
\beqnn
a(\eta_1^{-1})C_{i-1,1}+b(\eta_1^{-1})C_{i,1}+c(\eta_1^{-1})C_{i+1,1}=0,\quad i\geq1.
\eeqnn
Since $C_{i,1}$, $i\geq1$ satisfies the above relations, it should take the form of
\beqnn
C_{i,1}=A_1\bigg(\frac{1}{x_1(\eta_1^{-1})}\bigg)^{i-1}+B_1\bigg(\frac{1}{x_0(\eta_1^{-1})}\bigg)^{i-1},\quad i\geq1.
\eeqnn
To determine the coefficients $A_1$ and $B_1$, we solve the following initial equations:
\beqlb\label{5-12}
\left\{\begin{array}{ll}
b_2(\eta_1^{-1})C_{l,1}+\eta_1^{-1}c(\eta_1^{-1})(A_1+B_1)=0,
\\a(\eta_1^{-1})C_{l,1}+b(\eta_1^{-1})(A_1+B_1)+c(\eta_1^{-1})\bigg[A_1\bigg(\frac{1}{x_1(\eta_1^{-1})}\bigg)+B_1\bigg(\frac{1}{x_0(\eta_1^{-1})}\bigg)\bigg]=0.
\end{array}
\right.
\eeqlb
Let  $y=\eta_1^{-1}$, then
\beqlb\label{5-a1}
H_2(x_0(\eta_1^{-1}),\eta_1^{-1})=0.
\eeqlb
By  (\ref{2-5}) and (\ref{5-a1}), one can easily get that
\beqlb\label{5-a2}
h_2(x_0(\eta_1^{-1}),\eta_1^{-1})=\eta_1^{-1}a(\eta_1^{-1})x_0(\eta_1^{-1})+b_2(\eta_1^{-1})=0,
\eeqlb
since $h\big(x_0(y),y\big)=y\big[a(y)x_0(y)^2+b(y)x_0(y)+c(y)\big]=0$.
So, by (\ref{5-12}), (\ref{5-a2}) and  the equation that $a(y)x_0(y)x_1(y)=c(y)$, we get that
\beqnn
\left\{\begin{array}{ll}
A_1=-\frac{\eta_1b_2(\eta_1^{-1})}{c(\eta_1^{-1})}C_{l,1},
\\B_1=0.
\end{array}
\right.
\eeqnn
On the other hand, since $T(\frac{1}{\eta_1})=0$, we have
\beqlb\label{5-a3}
\sqrt{\Delta(\eta_1^{-1})}=\eta_1F(\eta_1^{-1}).
\eeqlb
According to the above arguments, we can get $x_1(\eta_1^{-1})$ and finally \eqref{5-9} by Theorem \ref{4-thm1}.

(ii) The proof for this case is similar to the case (i). Thus, it is omitted.

(iii) In this case, $\Delta{(y_0)}=0$ and then $x_0(y_0)=x_1(y_0)$. We assume that
$$\lim_{y\rightarrow y_0}\sqrt{1-\frac{y}{y_0}}\frac{d\psi_i(y)}{dy}=C_{i,3},\quad i\geq1.$$
By using the same method as in the case (i), we get
\beqnn
b_2(y_0)C_{l,3}+y_{0}c(y_0)C_{1,3}=0,
\eeqnn
and
\beqnn
a(y_0)C_{i-1,3}+b(y_0)C_{i,3}+c(y_0)C_{i+1,3}=0,\quad i\geq1.
\eeqnn
The solution is
\beqnn
C_{i+1,3}=(A_3+B_{3}i)\bigg(\frac{1}{x_1(y_0)}\bigg)^i,\quad i\geq0.
\eeqnn
By using the same method as in the case (i), we can determine $A_3$, $B_3$ and $x_1(y_0)$ and finally get  $C_{i,3}$.
%\beqnn
%A_3&=&-\frac{b_2(y_0)}{y_{0}c(y_0)}C_{l,3},
%\\B_3&=&-\frac{h_2(x_0(y_0),y_0)}{y_{0}c(y_0)}C_{l,3},
%\eeqnn
\qed

\subsection{Exact tail asymptotics for the marginal distribution}

In this subsection, we provide the details for the exact tail asymptotics of the marginal distribution $\pi_j^{(l)}=\sum_{i}\pi_{i,j}$, which can be characterized by computing $P(1,y)$.

From \eqref{2-10}, we have
\beqnn
P(1,y)=\frac{\bar{p}\mu_l}{qy}\big[(qy+\bar{q})P_2(y)-\bar{q}\pi_{00}\big].
\eeqnn
Hence
\beqnn
\pi_j^{(l)}=\bar{p}\mu_l\pi_{0,j},
\eeqnn
and the exact tail asymptotics for the marginal distribution $\pi_j^{(l)}$ is the same as that of the boundary probability $\pi_{0,j}$.

\section{Exact tail asymptotics for the high-priority queue}
In this section, we characterize the exact tail asymptotics of the joint stationary distribution as well as the marginal distribution along the high-priority queue direction. The former, $\pi_{i,j}$ for fixed $j\geq0$, can be derived from the generating function of $\varphi_j(x)$ given by \eqref{3-12}, and the latter $\sum_j{\pi_{i,j}}$ can be obtained from the fundamental form \eqref{2-10}.

\begin{lem}\label{l-6}
For $j\geq0$,
\beqlb\label{6-2}
\varphi_j(x)\sim C^j\pi_{0,0}(1-r_{0}x)^{-j-1},\quad as~ r_{0}x\rightarrow1,
\eeqlb
where \beqlb\label{6-1}
C=\frac{q}{\bar{q}}\frac{\big[\bar{p}\mu_{h}r_{0}^2+(p\mu_h+\bar{p}\bar{\mu}_h)r_0+p\bar{\mu}_h\big]}{p\bar{\mu}_h-\bar{p}\mu_{h}r_{0}^2}.
\eeqlb
\end{lem}
{\it Proof:} If $j=0$, then by \eqref{3-12}, one can easily get that the lemma holds.

Next, we will prove the lemma for $j\geq 1$ by induction.  By \eqref{3-12}, we have
\beqnn
\varphi_j(x)=\frac{a_j}{x-x_1}-\frac{q\varphi_{j-1}(x_0)(x+x_0)}{\bar{q}(x-x_1)}-\frac{q(px+\bar{p})(\bar{\mu}_{h}x+\mu_h)}{p\bar{q}\bar{\mu}_h(x-x_1)}\frac{\varphi_{j-1}(x)-\varphi_{j-1}(x_0)}{x-x_0},
\eeqnn
where $a_j$ and $\varphi_{j-1}(x_0)$ are constants depending on $j$.  It follows from  Lemma \ref{l-1} that for $j=1$,
\beqlb\label{3-a1}
\displaystyle\lim_{x\rightarrow x_1}\bigg(1-\frac{x}{x_1}\bigg)^2\varphi_1(x)=\frac{q(px_1+\bar{p})(\bar{\mu}_{h}x_1+\mu_h)}{p\bar{q}\bar{\mu}_h(x_1-x_0)x_1}\pi_{0,0}=C\pi_{0,0}.
\eeqlb
Substituting $x_0=\frac{r_0}{w}$, $x_1=\frac{1}{r_0}$ and $w=\frac{p\bar{\mu}_h}{\bar{p}\mu_{h}}$ into the  equation \eqref{3-a1} and after some elementary manipulations, we get \eqref{6-1}.

Now, we assume that \eqref{6-2} is true for $j=k$, then by induction
\beqnn
\displaystyle\lim_{x\rightarrow x_1}\bigg(1-\frac{x}{x_1}\bigg)^{j+2}\varphi_{j+1}(x)&=&-\frac{q(px_1+\bar{p})(\bar{\mu}_{h}x_1+\mu_h)}{p\bar{q}\bar{\mu}_h(x_1-x_0)x_1}\displaystyle\lim_{x\rightarrow x_1}\frac{\varphi_{j}(x)}{x-x_1}(x_1-x)^{j+2}\frac{1}{x_{1}^{j+1}}
\\&=&C\displaystyle\lim_{x\rightarrow x_1}\varphi_{j}(x)\big(1-\frac{x}{x_1}\big)^{j+1}
\\&=&C^{j+1}\pi_{0,0}.
\eeqnn
This completes the proof.
\qed

\begin{thm}
For the discrete-time preemptive priority queue with two classes of customers satisfying $\rho<1$, the exact tail asymptotics in the joint stationary distribution along the high-priority queue direction is characterized as follows: for a fixed number $j$ of low-priority customers, where $j\geq0$,
\beqnn
\pi_{i,j}\sim \frac{C^j\pi_{0,0}}{j!}i^{j}r_0^{i},
\eeqnn
where
$r_0$, $\pi_{0,0}$ and $C$ are given by \eqref{3-10}, \eqref{3-13} and \eqref{6-1}, respectively.
\end{thm}
{\it Proof:}
It is obvious that for $j\geq0$, $\varphi_j(x)$ is analytic in the region $\Delta(\phi,\epsilon)=\big\{x:|r_{0}x|<1+\epsilon, r_{0}x\neq 1,|Arg(r_{0}x-1)|>\phi,\epsilon>0,0<\phi<\pi/2,\big\}$. By Theorem \ref{4-thm1} and Lemma \ref{l-6},  we have for fixed $j$, where $j\geq0$,
\beqnn
\pi_{i,j}\sim\frac{C^{j}\pi_{0,0}}{\Gamma(j+1)}i^{j}\bigg(\frac{1}{r_{0}}\bigg)^{-i}=\frac{C^j\pi_{0,0}}{j!}i^{j}r_{0}^i.
\eeqnn
\qed

Finally, we present the tail asymptotic property of the marginal distribution $\pi_i^{(h)}=\sum_{j}\pi_{i,j}$. Taking $y=1$ in the fundamental form \eqref{2-10}, we have
\beqlb\label{6-a2}
P(x,1)=\frac{\bar{p}\mu_h\psi_0(1)}{\bar{p}\mu_h-p\bar{\mu}_{h}x}.
\eeqlb
On the other hand, by \eqref{3-20},
\beqlb\label{6-a1}
\psi_0(1)=\frac{1-\rho}{\bar{p}}.
\eeqlb
By \eqref{6-a2} and \eqref{6-a1},
\beqnn%\label{6-a3}
P(x,1)=\frac{\mu_h(1-\rho)}{\bar{p}\mu_h-p\bar{\mu}_{h}x},
\eeqnn
and from the assumption of this model, we have
\beqnn%\label{6-a4}
\frac{\bar{p}\mu_h}{p\bar{\mu}_h}>1.
\eeqnn
Therefore, $\frac{\bar{p}\mu_h}{p\bar{\mu}_h}$ is the dominant singularity and
\beqnn
\pi_{i}^{(h)}\sim\frac{1-\rho}{\bar{p}}\bigg(\frac{p\bar{\mu}_h}{\bar{p}\mu_h}\bigg)^{i}.
\eeqnn

\vspace*{5mm}
\noindent{\bf Acknowledgments:} The authors would like to thank the  referee for the helpful suggestions and comments which have improved the paper. This work was done during our visit
to the School of Mathematics and Statistics of Carleton University (Ottawa, Canada). The authors thank the school for providing us with a good working condition during our visit. The authors would also like to thank Professor Yiqiang Q. Zhao, Carleton University, for stimulating discussions. This work was supported by the National Natural Science Foundation of China (11271373, 11361007) and the Guangxi Natural Science Foundation (No.2012GXNSFBA053010,2014GXNSFCA118001). Song, Y. and Liu, Z. would also like to thank China Scholarship Council for supporting their visit to Carleton University.

\end{document}